\documentclass[a4paper,11pt]{article}
\usepackage[english]{babel}
\usepackage[latin1]{inputenc}
\usepackage{amssymb}
\usepackage{amsmath}
\usepackage{amscd}
\usepackage{latexsym}
\usepackage{amsthm}
\usepackage{eucal}
\usepackage{mathrsfs}
\theoremstyle{plain}
\newtheorem{thm}{Theorem}[section]
\newtheorem{cor}[thm]{Corollary}
\newtheorem{lem}[thm]{Lemma}
\newtheorem{prop}[thm]{Proposition}

\theoremstyle{definition}

\theoremstyle{remark}
\newtheorem{oss}{Remark}[section]
\newtheorem{esem}{Example}[section]

\begin{document}

\linespread{1.1}

\title{The $2-$Factoriality of the O'Grady
Moduli Spaces}

\author{Arvid Perego}

\maketitle

\begin{abstract}
The aim of this work is to show that the moduli space $M_{10}$
introduced by O'Grady in \cite{OG1} is a $2-$factorial variety.
Namely, $M_{10}$ is the moduli space of semistable sheaves with
Mukai vector $v:=(2,0,-2)\in H^{ev}(X,\mathbb{Z})$ on a projective
K3 surface $X$. As a corollary to our construction, we show that the 
Donaldson morphism gives a Hodge isometry between $v^{\perp}$ 
(sublattice of the Mukai lattice of $X$) and its image in $H^{2}
(\widetilde{M}_{10},\mathbb{Z})$, lattice with respect to the Beauville 
form of the $10-$dimensional irreducible symplectic manifold 
$\widetilde{M}_{10}$, obtained as symplectic resolution of $M_{10}$. 
Similar results are shown for the moduli space $M_{6}$ introduced by 
O'Grady in \cite{OG2}.
\end{abstract}

\section{Introduction}

Moduli spaces of semistable sheaves on abelian or projective K3 surfaces
are one of the main tools to produce examples of irreducible
symplectic manifolds. If $M_{v}$ denotes the moduli space of semistable
sheaves with Mukai vector $v$ on a projective K3 surface, 
it is a well-known result that if $v$ is primitive and the chosen polarization 
is $v-$generic, then $M_{v}$ is an irreducible symplectic manifold. Moreover,
$M_{v}$ is deformation equivalent to an Hilbert scheme of points on some 
projective K3 surface. An analogous result shows that if the surface is abelian, 
from $M_{v}$ one can produce an irreducible symplectic manifold, which is 
deformation equivalent to a generalized Kummer variety on some abelian surface.

The choice of a non-primitive Mukai vector can give rise to new examples.
Let $X$ be a projective K3 surface, and suppose there is an ample
divisor $H$ on $X$ such that $H^{2}=2$ and $Pic(X)=\mathbb{Z}\cdot
H$. Let us consider the moduli space $M_{10}$ of $H-$semistable sheaves
on $X$ whose Mukai vector is $(2,0,-2)\in H^{2*}(X,\mathbb{Z})$.
The moduli space $M_{10}$ was introduced by O'Grady in \cite{OG1},
where he shows that $M_{10}$ admits a symplectic resolution 
$\widetilde{M}_{10}$, which is an irreducible symplectic manifold of 
dimension 10. In \cite{OG2}, O'Grady introduced a $6-$dimensional
moduli space $M_{6}$ of semistable sheaves on an abelian surface, showing
that it admits a symplectic resolution $\widetilde{M}_{6}$, which is an
irreducible symplectic manifold of dimension 6. In both cases, the obtained
manifold is not deformation equivalent to any other previously known example 
of irreducible symplectic manifold.

A natural question is if there are other moduli spaces
of semistable sheaves admitting a symplectic resolution and giving rise to
new irreducible symplectic manifolds. In \cite{L-S} and \cite{K-L-S}, 
the authors answered to the question: in
particular, in \cite{L-S} it is shown that if $v=2w$, where $w$ is
a primitive Mukai vector such that $(w,w)=2$, then $M_{v}$ admits
a symplectic resolution, obtained as the blow-up of $M_{v}$ along
its reduced singular locus. In \cite{K-L-S} it is
shown that if $v=mw$ for $m\in\mathbb{N}$ and $w$ a primitive
Mukai vector, such that $m>2$ or $m=2$ and $(w,w)>2$, then $M_{v}$ does not
admit any symplectic resolution. In this case, $M_{v}$ is a
locally factorial scheme.

The aim of this work is to describe singularities of $M_{10}$ and
$M_{6}$. Namely, we show the following:

\begin{thm}
\label{thm:m102fatt}The moduli spaces $M_{10}$ and $M_{6}$
are $2-$factorial projective varieties.
\end{thm}

The main ingredient in the proof of Theorem \ref{thm:m102fatt} is
the Le Potier morphism, which to certain classes in $K_{top}(X)$
associates a line bundle on a moduli space.
In the case of $M_{10}$, using results of \cite{R1} we show that the 
Weil divisor $B$ parameterizing non-locally free sheaves is not a
Cartier divisor. Le Potier's construction allows us to show that
$2B$ is Cartier. We will deduce the $2-$factoriality of $M_{10}$
from this.
The same ideas are used in the proof of the $2-$factoriality
of $M_{6}$, but here the problem is subtler: the exceptional divisor 
$\widetilde{\Sigma}$ of the symplectic resolution $\widetilde{M}_{6}$ 
is divisible by $2$, while this is not the case for $M_{10}$. 
This property implies the existence of a Weil divisor $D$ on $M_{6}$
such that $2D=0$. Using results of \cite{R2}, we show that
$D$ is not a Cartier divisor and that $M_{6}$ is in fact
$2-$factorial.

As a corollary to our construction, we show the following

\begin{thm}
\label{thm:isodon}Let $X$ be a projective K3 surface such that $Pic(X)=\mathbb{Z}\cdot H$
for some ample line bundle such that $H^{2}=2$, and let $v=(2,0,-2)\in\widetilde{H}(X,\mathbb{Z})$. 
Let $v^{\perp}\subseteq\widetilde{H}(X,\mathbb{Z})$ be the orthogonal to $v$ with 
respect to the Mukai pairing. There is a Hodge injective morphism 
$$f:v^{\perp}\longrightarrow H^{2}(\widetilde{M}_{10},\mathbb{Z})$$which gives an isometry
between $v^{\perp}$ (lattice with respect to the Mukai form) and its image in $H^{2}(\widetilde{M}_{10},
\mathbb{Z})$ (lattice with respect to the Beauville form).
\end{thm}

An analogous result holds in the 6-dimensional example. 
This is the generalization of Theorem 0.1 in \cite{Y} for moduli spaces $M_{v}$, with $v$ primitive.

The paper is organized as follows. Sections from 2 to 5 are devoted
to the $10-$dimensional O'Grady's example: in section 2 we recall
the construction of $M_{10}$, and we show that it cannot be
locally factorial. In section 3 and 4 we show that $M_{10}$ 
is $2-$factorial, and in section 5 we prove Theorem \ref{thm:isodon}.

Sections from 6 to 9 are devoted to the $6-$dimensional O'Grady's
example, following the same structure described for the
previous example. The exposition of the two examples is presented as
symmetric as possible, and the main proofs for the $6-$dimensional case
are almost identical to those of the $10-$dimensional one. Anyway,
subtle differencies are shown when necessary.
Finally, in section 10 we present a brief appendix on constructions
of flat families that we need all along the paper.

\section{The local factoriality of $M_{10}$}

In this section we recall the construction of the $10-$dimensional 
moduli space $M_{10}$ and its main properties, namely
those contained in \cite{R1}. We provide two construction of
flat families of sheaves that we will use in sections 3 and 4, and
we show that $M_{10}$ is not locally factorial.

\subsection{Generalities on $M_{10}$}

Let us recall the setting of \cite{OG1}. Let $X$ be a projective
K3 surface such that $Pic(X)=\mathbb{Z}\cdot H$, where $H$ is an
ample line bundle with $H^{2}=2$. Let $M_{10}$ be the
moduli space of $H-$semistable sheaves on $X$ with Mukai vector
$(2,0,-2)$. It is a $10-$dimensional projective variety whose regular 
locus is $M_{10}^{s}$, the open subset parameterizing stable sheaves. 
Let $\Sigma$ be the singular locus of $M_{10}$, which is a codimension 
2 closed subset in $M_{10}$ (see \cite{OG1}). As semistable locally 
free sheaves are stable in this setting (by Lemma 1.1.5 in \cite{OG1}), 
the open subset $M_{10}^{lf}$ of $M_{10}$ parameterizing locally free
sheaves is contained in $M_{10}^{s}$. Let $B$ be the closed subset 
of $M_{10}$ parameterizing non locally free sheaves: then $\Sigma
\subseteq B$, and by \cite{OG1}, $B$ is an irreducible 
Weil divisor. The first result we need is:

\begin{thm}
\label{thm:ogls}(\textbf{O'Grady}, \textbf{Lehn-Sorger})
The moduli space $M_{10}$ admits a symplectic resolution
$\pi:\widetilde{M}_{10}\longrightarrow M_{10},$ and $\widetilde
{M}_{10}$ is a $10-$dimensional irreducible symplectic manifold
with $b_{2}\geq 24$. Moreover,
$\widetilde{M}_{10}$ can be obtained as the blow-up of $M_{10}$
along $\Sigma$ with its reduced schematic structure.
\end{thm} 

\proof The proof is in \cite{OG1}. In \cite{L-S} it is 
proved that $\widetilde{M}_{10}$ can be obtained as the blow up 
of $M_{10}$ along its reduced singular locus.\endproof

Let $\widetilde{\Sigma}$ be the exceptional divisor of $\pi$, and 
let $\widetilde{B}$ be the proper transform of $B$ under $\pi$. Let
$M_{10}^{\mu-ss}$ be the Donaldson-Uhlenbeck compactification of the
moduli space $M_{10}^{\mu}$ of $\mu-$stable sheaves, and
let $\phi:M_{10}\longrightarrow M_{10}^{\mu-ss}$ be the canonical surjective
morphism. As shown in \cite{OG1}, $M_{10}^{\mu-ss}=M_{10}^{lf}\coprod
S^{4}(X)$. Let $\delta$ be the fiber of $\pi$ over
a generic point in the smooth locus of $\Sigma$, and let $\gamma'$ 
be the fiber of $\phi$ over a generic point of the smooth locus of $S^{4}(X)\subseteq
M_{10}^{\mu-ss}$. Moreover, let $\gamma$ be the proper transform of $\gamma'$.
Finally, let $$\mu_{D}:H^{2}(X,\mathbb{Z})\longrightarrow
H^{2}(M_{10}^{\mu-ss},\mathbb{Z})$$be the Donaldson morphism (see \cite
{OG1}, \cite{F-M}).

\begin{thm}
\label{thm:rapagnetta2}(\textbf{Rapagnetta}). The second Betti
number of $\widetilde{M}_{10}$ is $24$. Moreover
\begin{enumerate}
\item The morphism
$\widetilde{\mu}:=\pi^{*}\circ\phi^{*}\circ\mu_{D}:H^{2}(X,\mathbb{Z})
\longrightarrow H^{2}(\widetilde{M}_{10},\mathbb{Z})$ is
injective. \item We have the following equalities:
$$c_{1}(\widetilde{\Sigma})\cdot\delta=-2,\,\,\,\,\,\,\,\,\,\,\,
c_{1}(\widetilde{B})\cdot\delta=1$$$$c_{1}(\widetilde{\Sigma})\cdot
\gamma=3,\,\,\,\,\,\,\,\,\,\,\,c_{1}(\widetilde{B})\cdot\gamma=-2.$$
\item The second integral cohomology of $\widetilde{M}_{10}$ is
$$H^{2}(\widetilde{M}_{10},\mathbb{Z})=\widetilde{\mu}(H^{2}(X,\mathbb{Z}))
\oplus\mathbb{Z}\cdot
c_{1}(\widetilde{\Sigma})\oplus\mathbb{Z}\cdot
c_{1}(\widetilde{B}).$$\item Let $q$ be the Beauville form of 
$\widetilde{M}_{10}$. Then for every $\alpha,\beta\in H^{2}(X,\mathbb{Z})$ we have
$$q(\widetilde{\mu}(\alpha),\widetilde{\mu}(\beta))=\alpha\cdot\beta,\,\,\,\,\,\,\,
q(\widetilde{\mu}(\alpha),c_{1}(\widetilde{\Sigma}))=q(\widetilde{\mu}(\alpha),
c_{1}(\widetilde{B}))=0,$$
$$q(c_{1}(\widetilde{\Sigma}),
c_{1}(\widetilde{\Sigma}))=-6,\,\,\,\,\,\,\,
q(c_{1}(\widetilde{\Sigma}),c_{1}(\widetilde{B}))=3,$$
$$q(c_{1}(\widetilde{B}),
c_{1}(\widetilde{\Sigma}))=3,\,\,\,\,\,\,\,
q(c_{1}(\widetilde{B}),c_{1}(\widetilde{B}))=-2.$$
\end{enumerate}
\end{thm}

\proof The proof is contained in \cite{R1}, Theorems 1.1, 3.1
and 4.3.\endproof

\subsection{Flat families}

In this subsection we present two examples of flat families of sheaves
we will use in the following. We refer to section 10 for the general
construction.

\begin{esem}
\label{esem:univx}Let $X$ be a projective K3 surface with
$Pic(X)=\mathbb{Z}\cdot H$, where $H$ is an ample line bundle such
that $H^{2}=2$. Fix three different points $x_{1},x_{2},x_{3}\in
X$, and consider
$$i:X\longrightarrow
S^{4}(X),\,\,\,\,\,\,i(x):=x+x_{1}+x_{2}+x_{3},$$which is a closed
immersion. Let $T:=i(X)\simeq X$, and consider a surjective morphism
$\varphi:\mathscr{O}_{X}^{2}\longrightarrow
\mathbb{C}_{x_{1}}\oplus\mathbb{C}_{x_{2}}\oplus\mathbb{C}_{x_{3}}$ as
in Proposition 3.0.5 in \cite{OG1}. Let
$\mathscr{K}:=\ker(\varphi),$ which is a rank 2 sheaf such that
$det(\mathscr{K})=\mathscr{O}_{X}$ and $c_{2}(\mathscr{K})=3$. 

Let $\Delta\subseteq T\times X$ be the diagonal (up to the isomorphism
between $T$ and $X$). By Corollary 2, 
Chapter II.5 in \cite{M} (see Lemma \ref{lem:mum} below), the sheaf
$p_{T*}\mathscr{H}om(p_{X}^{*} \mathscr{K},\mathscr{O}_{\Delta})$
is a rank 2 vector bundle, and for every $x\in T$ the canonical
morphism
$$(p_{T*}\mathscr{H}om(p_{X}^{*}
\mathscr{K},\mathscr{O}_{\Delta}))_{x}\longrightarrow Hom(
\mathscr{K},\mathbb{C}_{x})$$is an isomorphism. Let
$Y:=\mathbb{P}(p_{T*}\mathscr{H}om(p_{X}^{*}\mathscr{K},
\mathscr{O}_{\Delta}))$ and $p:Y\longrightarrow T$ be the canonical
projection. We have a canonical morphism (see Section 10)
$$\widetilde{f}:q_{X}^{*}\mathscr{K}\otimes
q_{Y}^{*}\mathscr{T}\longrightarrow (p\times
id_{X})^{*}\mathscr{O}_{\Delta},$$where $q_{X}$ and $q_{Y}$ are the natural
projections of $Y\times X$ to $X$ and $Y$ respectively, and $\mathscr{T}$ is the
tautological line bundle on $Y$. Consider $\mathscr{H}:=\ker(\widetilde{f})$.
\end{esem}

\begin{lem}
\label{lem:ssu}Let $\mathscr{E}$ be a sheaf defining a point in
$B$ and whose singular locus is given by $x,x_{1},x_{2},x_{3}$.
Then $\mathscr{E}$ defines a point $[f_{\mathscr{E}}]\in Y$, and the
restriction $\mathscr{H}_{[f_{\mathscr{E}}]}:=\mathscr{H}_{|q_{Y}^{-1}
([f_{\mathscr{E}}])}$ is isomorphic to $\mathscr{E}.$
Moreover, the morphism $\widetilde{f}$ is surjective and $\mathscr{H}$ is a
$Y-$flat family.
\end{lem}

\proof The sheaf $\mathscr{E}$ is the kernel of a surjective
(hence non-zero) morphism $f_{\mathscr{E}}:\mathscr{K}\longrightarrow\mathbb{C}_{x}$
(see section 3.1 in \cite{OG1}), defining a point
$t=[f_{\mathscr{E}}]\in p^{-1}(x)$ since
$p^{-1}(x)\simeq\mathbb{P}(Hom(\mathscr{K},\mathbb{C}_{x}))$. By
definition of $\widetilde{f}$, we have
$\widetilde{f}_{t}=f_{\mathscr{E}}$. The morphism $\widetilde{f}$
is surjective: indeed, co$\ker(\widetilde{f})$ is trivial if and
only if it is trivial on the fibers of $q_{Y}$. If $t\in Y$,
then $t$ corresponds to a surjective morphism $f_{\mathscr{E}}$,
so that co$\ker(\widetilde{f})_{t}=$ co$\ker( f_{\mathscr{E}})=0$,
and we are done.

Since $\widetilde{f}$ is surjective, the family $\mathscr{H}$ is
$Y-$flat. Now, since $q_{X}^{*}\mathscr{K}\otimes
q_{Y}^{*}\mathscr{T}$ and $(p\times
id_{X})^{*}\mathscr{O}_{\Delta}$ are $Y-$flat, for every $t\in Y$
the canonical morphism
$$\mathscr{H}_{t}\longrightarrow(q_{X}^{*}\mathscr{K}
\otimes q_{Y}^{*}\mathscr{T})_{t}\simeq\mathscr{K}$$is injective,
so that $\mathscr{H}_{t}=\ker(\widetilde{f}_{t})$. As $t\in Y$ corresponds to a
surjective morphism $f_{\mathscr{E}}:\mathscr{K}
\longrightarrow\mathbb{C}_{x}$, where $x=p(t)$, whose kernel is
$\mathscr{E}$, and $\widetilde{f}_{t}= f_{\mathscr{E}}$, we are done.
\endproof

\begin{esem}
\label{esem:piuno}In the same setting of the previous example, let
$x\in X$ be different from $x_{1},x_{2},x_{3}$, and let
 $T:=\{x\}$. Moreover, let $i$ be the
inclusion of $T$ in $X$. Let $Y:=\mathbb{P}(p_{T*}\mathscr{H}om
(p_{X}^{*}\mathscr{K},i_{*}\mathbb{C}_{x}))\simeq\mathbb{P}^{1}$,
By the general construction, we have
$\widetilde{f}:q_{X}^{*}\mathscr{K}\longrightarrow (p\times
id_{X})^{*}i_{*}\mathbb{C}_{x}\otimes q_{Y}^{*}\mathscr{O}(1),$
where $p:Y\longrightarrow T$ is the canonical morphism.
In particular, notice that $(p\times
id_{X})^{*}i_{*}\mathbb{C}_{x}\otimes q_{Y}^{*}\mathscr{O}(1)
=j_{*}\mathscr{O}(1),$ where
$j:\mathbb{P}^{1}\times T\longrightarrow \mathbb{P}^{1}\times
X$ is the inclusion. In conclusion, we have
$\widetilde{f}:q_{X}^{*}\mathscr{K}\longrightarrow j_{*}\mathscr{O}(1).$ 
Finally, let $\mathscr{H}:=\ker(\widetilde{f})$.
\end{esem}

\begin{lem}
\label{lem:piunos}Let $\mathscr{E}$ be a sheaf defining a point in
$B$ whose singular locus is given by $x,x_{1},x_{2},x_{3}$, and
let $[f_{\mathscr{E}}]$ be the point of $Y$ defined by
$\mathscr{E}$. Then
$\mathscr{H}_{[f_{\mathscr{E}}]}\simeq\mathscr{E},$ and
$\widetilde{f}$ is a surjective morphism. Moreover, the family
$\mathscr{H}$ is $Y-$flat.
\end{lem}

\proof The proof works as the one of Lemma \ref{lem:ssu}.\endproof

\subsection{The moduli space $M_{10}$ is not locally factorial}

A first application of Theorem \ref{thm:rapagnetta2} is the
following:

\begin{prop}
\label{prop:nolocfact}If $n\in\mathbb{Z}$ is such that $nB$
is a Cartier divisor, then $n$ is even. In particular,
$M_{10}$ is not locally factorial.
\end{prop}

\proof Let $n\in\mathbb{Z}$ be such that $nB$ is Cartier. Then
$\pi^{*}(nB)=n\widetilde{B}+m\widetilde{\Sigma}$ for some $m\in\mathbb{Z}$,
since $\widetilde{B}$ is the proper transform of $B$. By the
projection formula we have $c_{1}(\pi^{*}(nB))\cdot\delta=0$, as
$\delta$ is contracted by $\pi$. By point 2 of Theorem
\ref{thm:rapagnetta2}, we get
$$0=c_{1}(\pi^{*}(nB))\cdot\delta=nc_{1}(\widetilde{B})\cdot\delta
+mc_{1}(\widetilde{\Sigma})\cdot\delta=n-2m,$$so that $n$ is
even. Finally, $M_{10}$ is not locally factorial: if it was, 
then $B$ would be a Cartier divisor, which
is clearly not the case by the previous part of the proof.
\endproof

\begin{oss}
\label{oss:notor10}Theorem \ref{thm:rapagnetta2}
implies even that $Pic(M_{10})$ is free. Indeed, let $L\in
Pic(M_{10})$ be torsion of period $t\in\mathbb{N}$, and let
$\widetilde{L}$ be its proper transform under $\pi$. Then
$\pi^{*}(L)=\widetilde{L}+m\widetilde{\Sigma}$ for some
$m\in\mathbb{Z}$, and $t(\widetilde{L}+m\widetilde{\Sigma})=0.$ As
$\widetilde{M}_{10}$ is simply connected, by
point 3 of Theorem \ref{thm:rapagnetta2} we see that
$Pic(\widetilde{M}_{10})$ is free: in conclusion
$\widetilde{L}=-m\widetilde{\Sigma}$, so that $L=0$. The same
argument shows that $\pi^{*}:Pic(M_{10})\longrightarrow
Pic(\widetilde{M}_{10})$ is injective.

Moreover, $c_{1}:Pic(M_{10})\longrightarrow
H^{2}(M_{10},\mathbb{Z})$ is injective: if $L,L'\in Pic(M_{10})$ are
such that $c_{1}(L)=c_{1}(L')$, then $c_{1}(\pi^{*}(L))=c_{1}(\pi^{*}(L'))$,
so that $\pi^{*}(L)=\pi^{*}(L')$. As $\pi^{*}$ is injective on $Pic(M_{10})$, this
implies $L\simeq L'$.
\end{oss}

To conclude this section, we show the following:

\begin{lem}
\label{lem:intgamma}If $n\in\mathbb{Z}$ is such that $nB$ is
a Cartier divisor, then
$$c_{1}(nB)\cdot\gamma'=-\frac{n}{2}\in\mathbb{Z}.$$
\end{lem}

\proof The fact that $-n/2$ is an integer follows from
Proposition \ref{prop:nolocfact}, as $nB$ is a Cartier divisor.
By definition of $\gamma$ and $\gamma'$, there is $l\in \mathbb{Q}$ 
such that $\pi^{*}(\gamma\,')=\gamma+ l\delta$. By point 2 in Theorem
\ref{thm:rapagnetta2} and the projection formula we have
$$3=c_{1}(\widetilde{\Sigma})\cdot\gamma=c_{1}(\widetilde{\Sigma})
\cdot\pi^{*}(\gamma\,')-l(c_{1}(\widetilde{\Sigma})\cdot\delta)=2l,$$
so that $\pi^{*}(\gamma\,')=\gamma+\frac{3}{2}\delta.$ Now, suppose that
$n\in\mathbb{Z}$ is such that $nB$ is a Cartier divisor. By
the projection formula $c_{1}(nB)\cdot\gamma'=
nc_{1}(\widetilde{B})\cdot\pi^{*}(\gamma'),$ so that
$$c_{1}(nB)\cdot\gamma'=nc_{1}(\widetilde{B})\cdot\gamma+
\frac{3n}{2}c_{1}(\widetilde{B})\cdot\delta=-2n+\frac{3n}{2}=-\frac{n}{2},$$by
point 2 of Theorem \ref{thm:rapagnetta2}, and we are done.
\endproof

\section{Line bundles on $M_{10}$}

In this section we study properties of line bundles on $M_{10}$. 
The main ingredients are Le Potier's construction
of line bundles on moduli spaces of sheaves on algebraic surfaces,
and the construction of flat families we presented in section 2.2.

\subsection{Le Potier's construction}

We recall Le Potier's construction (see \cite{LeP}
or \cite{H-L}, Chapter 8). Let $S$ be a
Noetherian scheme, and let
$\mathscr{F}$ be an $S-$flat family on $S\times X$. We can
define a morphism
$$\widetilde{\lambda}:K_{top}(X)\longrightarrow
Pic(S),\,\,\,\,\,\,\,\widetilde{\lambda}(\alpha):=det(p_{R!}(p_{X}^{*}\alpha\cdot
[\mathscr{F}]).$$We apply this construction when
$S$ is the open subset $R$ of a Grothendieck Quot-scheme whose quotient is 
$M_{10}$, and $\mathscr{F}$ is a universal family on $R\times X$.
Let $e:=[\mathscr{E}]\in K_{top}(X)$ be the class of a sheaf $\mathscr{E}$ 
parameterized by $M_{10}$, $h:=[H]\in K_{top}(X)$ and let $$\xi:K_{top}(X)\times
K_{top}(X)\longrightarrow\mathbb{Z},\,\,\,\,\,\,\,\,\xi(\alpha,\beta):=
\chi(\alpha\cdot\beta).$$By Theorem 8.1.5 in \cite{H-L},  
$\widetilde{\lambda}(\alpha)\in Pic(R)$ descends to a line bundle
$\lambda(\alpha)\in Pic(M_{10})$ if $\alpha\in
e^{\perp}\cap\{1,h,h^{2}\}^{\perp\perp}$ (orthogonality with respect
to $\xi$), so that there is a group morphism
$$\lambda:e^{\perp}\cap\{1,h,h^{2}\}^{\perp\perp}\longrightarrow
Pic(M_{10}).$$The first result we need is:
 
\begin{lem}
\label{lem:perp}Let $\alpha\in K_{top}(X)$. Then $\alpha\in e^{\perp}\cap
\{1,h,h^{2}\}^{\perp\perp}$ if and only if 
$c_{1}(\alpha)\in Pic(X)$ and $ch_{2}(\alpha)=0$.
\end{lem}

\proof By the Hirzebruch-Riemann-Roch Theorem, a class $\beta\in 
K_{top}(X)$ is in $\{1,h,h^{2}\}^{\perp}$ if and only if $v(\beta)= 
(0,b,0)$, where $b\in H^{2}(X,\mathbb{Z})$ is such that $b\cdot 
c_{1}(H)=0$. Then $\beta\in H^{2}(X,\mathbb{Z})\cap(H^{2,0}(X)\oplus H^{0,2}(X))$,
as $Pic(X)=\mathbb{Z}\cdot H$.

Now, let $\alpha\in K_{top}(X)$. Then
$\alpha\in\{1,h,h^{2}\}^{\perp\perp}$ if and only if
$\chi(\alpha\cdot\beta)=0$ for every
$\beta\in\{1,h,h^{2}\}^{\perp}$. By the previous part, we get
$c_{1}(\alpha)\cdot b=0$ for every $b\in H^{2}(X,\mathbb{Z})\cap
(H^{2,0}(X)\oplus H^{0,2}(X))$. Then $c_{1}(\alpha)$ has to be the first 
Chern class of a line bundle on $X$.
Finally, by the Hirzebruch-Riemann-Roch Theorem, we have $\alpha\in e^{\perp}$ if
and only if $ch_{2}(\alpha)=0$, as $ch(e)=(2,0,-4)$.
\endproof

Using this lemma, we are able to prove the following:

\begin{prop}
\label{prop:lineb}Let $p\in X$ be any point, and let $$u_{1}:Pic(X)\longrightarrow e^{\perp}\cap
\{1,h,h^{2}\}^{\perp\perp},\,\,\,\,\,\,\,\,u_{1}(L):=
[\mathscr{O}_{X}-L]+\frac{c_{1}^{2}(L)}{2}[\mathbb{C}_{p}],$$
$$u_{2}:\mathbb{Z}\longrightarrow
e^{\perp}\cap\{1,h,h^{2}\}^{\perp\perp},\,\,\,\,\,\,\,\,\, u_{2}(n):=
n[\mathscr{O}_{X}].$$Then $u:=u_{1}+u_{2}$ is a group isomorphism.
\end{prop}

\proof Let $L\in Pic(X)$. The Mukai vector of $u_{1}(L)$ 
is $(0,-c_{1}(L),0),$ so
that $u_{1}(L)\in
e^{\perp}\cap\{1,h,h^{2}\}^{\perp\perp}$ by Lemma \ref{lem:perp}. 
Moreover, for every $L_{1},L_{2}\in Pic(X)$
we have $v(u_{1}(L_{1}\otimes L_{2}))=v(u_{1}(L_{1})+u_{1}(L_{2})),$ 
where $v:K_{top}(X)\longrightarrow H^{2*}(X,\mathbb{Z})$ 
is the morphism sending a class in $K_{top}(X)$ to its Mukai vector.
As $v$ is a group isomorphism (see \cite{K}), then
$u_{1}$ is a group morphism.

Let $n\in\mathbb{Z}$. The Mukai vector of $u_{2}(n)$ is 
$v(u_{2}(n))=(n,0,n)$, so that $u_{2}(n)\in e^{\perp}\cap\{1,h,h^{2}\}^
{\perp\perp}$ by Lemma \ref{lem:perp}, and $u_{2}$ is clearly a group morphism.

We need to show that $u$ is an isomorphism. For the injectivity,
let $(L,n),(M,m)\in Pic(X) \oplus\mathbb{Z}$ be such that
$u(L,n)= u(M,m)$. Their Mukai
vectors are then equal: as these are, respectively, $(n,-c_{1}(L),n)$ and
$(m,-c_{1}(M),m)$, this implies $n=m$ and $c_{1}(L)=c_{1}(M)$. As
$X$ is a K3 surface, this implies $L\simeq M$, and
injectivity is shown. For the surjectivity, let $\alpha\in
e^{\perp}\cap\{1,h,h^{2}\}^{\perp\perp}$: by Lemma
\ref{lem:perp}, we have
$v(\alpha)=(r,c_{1}(L),r)$ for some $r\in\mathbb{Z}$ and $L\in
Pic(X)$. Then $v(\alpha)=v(u(L^{-1},r))$, and $\alpha=u(L^{-1},r)$.
\endproof

\begin{prop}
\label{prop:noint}We have the following intersecion properties. 
\begin{enumerate}
\item For every $L\in Pic(X)$, we have $c_{1}(\lambda(u_{1}(L)))
\cdot\gamma'=0$. \item For every $n\in\mathbb{Z}$ we have 
$c_{1}(\lambda(u_{2}(n)))\cdot\gamma'=-n$.
\end{enumerate}
\end{prop}

\proof We begin with the first item. As $Pic(X)=\mathbb{Z}\cdot H$, 
we need to verify the statement only for $H$. By Proposition 8.2.3 
in \cite{H-L} there is $m\in\mathbb{N}$ such that 
$\lambda(u(H))^{\otimes m}$ is generated by its global sections, 
and the canonical map
$$\phi:M_{10}\longrightarrow\mathbb{P}
(H^{0}(M_{10},\lambda(u(H))^{\otimes m})^{*})$$ has $M_{10}^{\mu
-ss}$ as image. In particular,
$\phi^{*}\mathscr{O}(1)=\lambda(u(H))^{\otimes m}$, so that
$$mc_{1}(\lambda(u(H)))\cdot\gamma'=c_{1}(\lambda(u(H))^{\otimes m})
\cdot\gamma'=c_{1}(\phi^{*}\mathscr{O}(1))\cdot\gamma'=0,$$as
$\gamma'$ is contracted by $\phi$. Finally
$c_{1}(\lambda(u(H)))\cdot\gamma'=0$, and we are done.

For the second item, we need to verify the statement only for
$n=1$. Notice that $c_{1}(\lambda(u_{2}(1)))\cdot\gamma'=
c_{1}(\lambda([\mathscr{O}_{X}])_{|\gamma'})$. Using the family
$\mathscr{H}$ defined in Example \ref{esem:piuno}, we have
$$c_{1}(\lambda([\mathscr{O}_{X}])_{|\gamma'})=c_{1}(q_{Y!}
(q_{X}^{*}[\mathscr{O}_{X}]\cdot[\mathscr{H}])).$$By the 
Grothendieck-Riemann-Roch Theorem, as the fibers of $q_{Y}$ are of
dimension 2 we have
$$c_{1}(q_{Y!}(q_{X}^{*}[\mathscr{O}_{X}]\cdot[\mathscr{H}]))=
q_{Y*}[q_{X}^{*}(ch(\mathscr{O}_{X})td(X)^{-1})\cdot
ch(\mathscr{H})]_{3}=$$ $$=-2q_{X}^{*}[y]\cdot
ch_{1}(\mathscr{H})+ch_{3}(\mathscr{H}),$$where $[y]$ is the class
of a point in $X$. By the Grothendieck-Riemann-Roch Theorem and
by definition of $\mathscr{H}$ we have
$ch_{1}(\mathscr{H})=[j_{*}ch(\mathscr{O}_{\mathbb{P}^{1}}(1))]_{1}$
and
$$ch_{3}(\mathscr{H})=[j_{*}ch(\mathscr{O}_{\mathbb{P}^{1}}(1))]_{3}
-2q_{X}^{*}[y]\cdot[j_{*}ch(\mathscr{O}_{\mathbb{P}^{1}}(1))]_{1}.$$In
conclusion, $ch_{1}(\mathscr{H})=0$ and
$ch_{3}(\mathscr{H})=-q_{Y}^{*}[p],$ where $[p]$ is the class of a
point in $Y$. Finally, we get
$$c_{1}(\lambda([\mathscr{O}_{X}]))\cdot\gamma'=q_{Y*}(-q_{Y}^{*}[p])=-1,$$
and we are done.\endproof

\subsection{Donaldson's and Le Potier's morphisms}

The aim of this section is to prove that the morphism $\lambda\circ u$
is injective. The main result we need is the following:

\begin{prop}
\label{prop:eqdonpot}Let $L\in Pic(X)$. Then
$c_{1}(\lambda(u_{1}(L)))=\phi^{*}\mu_{D}(c_{1}(L))$.
\end{prop}

\proof The proof is done in two steps: first we show that these
two classes are equal when restricted to a well-chosen subvariety;
then we show that the equality on this restriction implies the
equality everywhere.

\textit{Step 1}. Here, we refer to
Example \ref{esem:univx} for the notations. 
Consider the inclusion $i:T\simeq X\longrightarrow M_{10}^{\mu-ss}$
described in Example \ref{esem:univx}, and consider the morphism 
$j:Y\longrightarrow M_{10}$ induced by the family $\mathscr{H}$. By Lemma
\ref{lem:ssu}, $j$ is injective and its image is
$\phi^{-1}(T)$. For every $L\in Pic(X)$ we have
$$j^{*}\phi^{*}(\mu_{D}(c_{1}(L)))=p^{*}i^{*}(\mu_{D}(c_{1}(L))).$$By
Proposition 6.5 in \cite{F-M} we have
$i^{*}(\mu_{D}(c_{1}(L)))=c_{1}(L)\in NS(X)$ (up to the
isomorphism between $X$ and $T$), and we need to show that
$j^{*}c_{1}(\lambda(u_{1}(L)))=p^{*}(c_{1}(L))$. By Theorem 8.1.5 in
\cite{H-L} and Lemma \ref{lem:ssu}
$$j^{*}\lambda(u_{1}(L))=det(q_{Y!}(q_{X}^{*}u_{1}(L)\cdot[\mathscr{H}])),$$so
that by the Grothendieck-Riemann-Roch Theorem we get
$$c_{1}(j^{*}\lambda(u_{1}(L)))=
q_{Y*}[q_{X}^{*}(ch(u_{1}(L))td(X)^{-1})\cdot
ch(\mathscr{H})]_{3}=$$ $$=-q_{X}^{*}(c_{1}(L))\cdot
ch_{2}(\mathscr{H}).$$By Lemma \ref{lem:ssu} we finally have
$$q_{Y*}(-q_{X}^{*}(c_{1}(L))\cdot ch_{2}(\mathscr{H}))=q_{Y*}(q_{X}^{*}
(c_{1}(L))\cdot(p\times id_{X})^{*}[\Delta])=p^{*}(c_{1}(L)).$$

\textit{Step 2}. Let $L\in Pic(X)$ and
$\beta:=\phi^{*}\mu_{D}(c_{1}(L))-c_{1}(\lambda(u_{1}(L)))\in
H^{2}(M_{10},\mathbb{Z}).$ We need to show that $\beta=0$. 

By Step 1, $j^{*}\beta=0$. Moreover, $\beta\cdot\gamma'=0$: indeed,
$\phi^{*}\mu_{D}(c_{1}(L))\cdot\gamma'=0$ as $\gamma'$ is
contracted by $\phi$, and $c_{1}(\lambda(u_{1}(L)))\cdot\gamma'=0$ by
point 1 of Proposition \ref{prop:noint}. Now, by point 3 of Theorem
\ref{thm:rapagnetta2}, there are $\alpha\in H^{2}(X,\mathbb{Z})$
and $n,m\in\mathbb{Z}$ such that
$\pi^{*}\beta=\widetilde{\mu}(\alpha)+nc_{1}(\widetilde{\Sigma})
+mc_{1}(\widetilde{B}).$ By point 2 of Theorem
\ref{thm:rapagnetta2}, we get
$$0=\pi^{*}\beta\cdot\delta=\widetilde{\mu}(\alpha)\cdot\delta+
nc_{1}(\widetilde{\Sigma})\cdot\delta+mc_{1}(\widetilde{B})\cdot\delta=m-2n$$as
$\delta$ is contracted by $\pi$, and
$$0=\pi^{*}\beta\cdot\gamma=\widetilde{\mu}(\alpha)\cdot\gamma+
nc_{1}(\widetilde{\Sigma})\cdot\gamma+mc_{1}(\widetilde{B})\cdot\gamma=3n-2m$$since
$\pi^{*}\beta\cdot\gamma=\beta\cdot\gamma'=0$. In conclusion,
$n=m=0$ and $\pi^{*}\beta=\widetilde{\mu}(\alpha)$. This implies
$\beta=\phi^{*}\mu_{D}(\alpha)$: indeed, $\beta$ and
$\phi^{*}\mu_{D}(\alpha)$ are in $NS(M_{10})$, and $\pi^{*}$ is
injective on $NS(M_{10})$ by Remark \ref{oss:notor10}.
Restricting to $Y$ we then get
$$0=j^{*}\beta=j^{*}\phi^{*}\mu_{D}(\alpha)=p^{*}(\alpha),$$the last equality
coming from Proposition 6.5 in \cite{F-M}. To conclude, simply
note that $p^{*}:NS(T)\longrightarrow NS(Y)$ is injective as
$Y$ is a $\mathbb{P}^{1}-$bundle on $T$, so that 
$\alpha=0$, and we are done.
\endproof

\begin{cor}
\label{cor:pic}The morphism $\lambda\circ u
:Pic(X)\oplus\mathbb{Z}\longrightarrow Pic(M_{10})$ is injective.
Moreover, we have
$Pic(\widetilde{M}_{10})=\pi^{*}\circ\lambda\circ
u(Pic(X))\oplus\mathbb{Z}[\widetilde{\Sigma}]\oplus\mathbb{Z}
[\widetilde{B}].$
\end{cor}

\proof For the injectivity of $\lambda\circ u$, let
$L,M\in Pic(X)$ and $n,m\in\mathbb{Z}$ be such that
$\lambda(u(L,n))=
\lambda(u(M,m))$. By Proposition \ref{prop:noint} we have
$$-n=c_{1}(\lambda(u(L,n)))\cdot\gamma'=
c_{1}(\lambda(u(M,m)))\cdot\gamma'=-m,$$so that
$m=n$ and $\lambda(u_{1}(L))=\lambda(u_{1}(M))$. In particular,
$c_{1}(\lambda(u_{1}(L)))=c_{1}(\lambda(u_{1}(M)))$, so that $\phi^{*}
\mu_{D}(c_{1}(L))=\phi^{*}\mu_{D}(c_{1}(M))$ by Proposition \ref{prop:eqdonpot}.
Now, by point 1 of Theorem \ref{thm:rapagnetta2}, the morphism $\phi^{*}\circ
\mu_{D}$ is injective, so that we finally get $c_{1}(L)=c_{1}(M)$, implying
$L\simeq M$ as $X$ is a K3 surface.

By Remark \ref{oss:notor10} the morphism $\pi^{*}:Pic(M_{10})
\longrightarrow Pic(\widetilde{M}_{10})$ is injective. Moreover, as 
$\lambda\circ u$ is injective, the morphism $\lambda\circ u_{1}$ is, so that
$$\pi^{*}\circ\lambda\circ
u_{1}:Pic(X)\longrightarrow Pic(\widetilde{M}_{10})$$is injective.
To conclude, let $L\in Pic(\widetilde{M}_{10})$. By point 3 of Theorem
\ref{thm:rapagnetta2}, there are $\alpha\in H^{2}(X,\mathbb{Z})$
and $n,m\in\mathbb{Z}$ such that $c_{1}(L)=\widetilde{\mu}(\alpha)+nc_{1}(\widetilde{\Sigma})
+mc_{1}(\widetilde{B}).$ In particular $\widetilde{\mu}(\alpha)\in
NS(\widetilde{M}_{10})$, so that $\phi^{*}\mu_{D}(\alpha)\in NS(M_{10})$.
By Proposition 6.5 in \cite{F-M} $j^{*}\phi^{*}\mu_{D}(\alpha)=p^{*}(\alpha)$,
so that $p^{*}(\alpha)\in NS(Y)$. But this implies
$\alpha\in NS(X)$, and we are done.
\endproof

\section{The $2-$factoriality of $M_{10}$}

Using the results of the previous section, we are finally able to
show the $2-$factoriality of $M_{10}$.

\begin{prop}
\label{prop:picm}Let $A^{1}(M_{10})$ be the group of Weil divisors
of $M_{10}$ modulo linear equivalence. Then
$A^{1}(M_{10})=\lambda(u_{1}(Pic(X)))\oplus\mathbb{Z}[B].$
\end{prop}

\proof Notice that $A^{1}(M_{10})\simeq Pic(\pi^{-1}(M_{10}^{s})$. Indeed,
$\pi$ is an isomorphism on $M_{10}^{s}$, so that $Pic(\pi^{-1}(M^{s}_{10}))\simeq 
Pic(M_{10}^{s})$, and $Pic(M_{10}^{s})=A^{1}(M_{10}^{s})$, as $M_{10}^{s}$ is
smooth. Since $\Sigma=M_{10}\setminus M_{10}^{s}$ has codimension 2 in
$M_{10}$, we have $A^{1}(M_{10}^{s})=A^{1}(M_{10})$. 
Now, let us consider the sequence $$0\longrightarrow
\mathbb{Z}\stackrel{\sigma}\longrightarrow
Pic(\widetilde{M}_{10})\stackrel{\rho}\longrightarrow
Pic(\pi^{-1}(M_{10}^{s}))\longrightarrow 0,$$where $\sigma(1):=
[\widetilde{\Sigma}]$ and $\rho$ is the restriction morphism.
We claim that it is exact: $\sigma$ is clearly injective and
$\rho$ is surjective (see \cite{H}, Proposition 6.5). Moreover,
$\widetilde{\Sigma}\in\ker(\rho)$. Let $L\in\ker(\rho)$:
we need to show that $L$ is a multiple of
$\widetilde{\Sigma}$. By Corollary \ref{cor:pic} there are $M\in
Pic(X)$ and $n,m\in\mathbb{Z}$ such that
$$L=\pi^{*}(\lambda(u_{1}(M)))+n\widetilde{B}+m\widetilde{\Sigma}.$$As
$\rho(\widetilde{\Sigma})=0$, we get
$\rho(L)=\lambda(u_{1}(M))+nB$ as Weil divisors on
$M_{10}$, so that $nB=-\lambda(u_{1}(M))$. In particular, $nB$ is a 
Cartier divisor: by Lemma
\ref{lem:intgamma} and point 1 of Proposition \ref{prop:noint}, we have
$$-\frac{n}{2}=nB\cdot\gamma'= \lambda(u_{1}(M))\cdot\gamma'=0,$$so
that $n=0$ and $M=\mathscr{O}_{X}$ (by Proposition
\ref{prop:eqdonpot}). By Corollary \ref{cor:pic} we are done.
\endproof

\begin{cor}
\label{cor:bnur}The only Weil divisors on $M_{10}$ that are
possibly not Cartier are the multiples of $B$.
\end{cor}

\proof As $\lambda(u_{1}(Pic(X)))\subseteq Pic(M_{10})$, 
this is an immediate corollary of Proposition \ref{prop:picm}. 
\endproof

The final result of this section is the following, which is one of the
main results of the paper:

\begin{thm}
\label{thm:2b}There is $L\in Pic(X)$ such that $2B=\lambda
(u(L,1))\in Pic(M_{10}).$ In
particular, the moduli space $M_{10}$ is $2-$factorial.
\end{thm}

\proof As $B$ is not a Cartier divisor by Proposition \ref{prop:nolocfact},
by Corollary \ref{cor:bnur} the $2-$factoriality of $M_{10}$ follows
once we show $2B\in Pic(M_{10})$. 
By Proposition \ref{prop:picm}, there are $n\in\mathbb{Z}$ and $M\in Pic(X)$ such that
$\lambda(u_{2}(1))=\lambda(u_{1}(M))+nB.$ In particular $nB\in Pic(M_{10})$, 
so that by Proposition \ref{prop:noint} and Lemma \ref{lem:intgamma}
$$-1=c_{1}(\lambda(u_{2}(1)))\cdot\gamma'=
c_{1}(\lambda(u_{1}(M)))\cdot\gamma'+c_{1}(nB)\cdot\gamma'=-\frac{n}{2}.$$In 
conclusion $n=2$ and $2B=\lambda(u(M^{-1},1))$.
\endproof

\section{The Beauville form of $\widetilde{M}_{10}$}

The aim of this section is to show that the line bundle $L$ in the
statement of Theorem \ref{thm:2b} is trivial. As a consequence of this, we
prove Theorem \ref{thm:isodon}.

\subsection{Properties of the Weil divisor $B$}

In this section we show some properties of the sheaves
parameterized by $M_{10}$. The main result is that the Weil
divisor $B$ is characterized in cohomological terms. We begin with
two lemmas.

\begin{lem}
\label{lem:edual}Let $E$ be a rank 2 locally free sheaf with
trivial determinant. Then $E\simeq E^{*}$.
\end{lem}

\proof By hypothesis on $E$, the canonical
morphism $E\otimes E\longrightarrow E\wedge E$ is a perfect
pairing.
\endproof

\begin{lem}
\label{lem:no0}Let $\mathscr{E}$ be a sheaf defining a point in
$M_{10}$. Then $H^{0}(X,\mathscr{E})=0$ and
$h^{1}(X,\mathscr{E})=h^{2}(X,\mathscr{E})$.
\end{lem}

\proof By the Hirzebruch-Riemann-Roch Theorem, the Hilbert
polynomial of $\mathscr{E}$ is $P(\mathscr{E},n)=2n^{2}$, since
$H^{2}=2$. In particular, $\chi(\mathscr{E})=0$, so that
$h^{1}(X,\mathscr{E})=h^{0}(X,\mathscr{E})+h^{2}(X,\mathscr{E})$.
We show that $h^{0}(X,\mathscr{E})=0$. Recall that the 
reduced Hilbert polynomial of a sheaf $\mathscr{F}$ of
dimension 2 on a surface is $p(\mathscr{F},n):=P(\mathscr{F},n)/rk(\mathscr{F})$:
then $p(\mathscr{E},n)=n^{2}$, and
$p(\mathscr{O}_{X},n)=n^{2}+2$. By Proposition 1.2.7 in
\cite{H-L}, we have $Hom(\mathscr{O}_{X},\mathscr{E})=0$, 
and we are done.\endproof

As a consequence, we have the following:

\begin{prop}
\label{prop:coholf}Let $\mathscr{E}$ be a semistable sheaf 
with Mukai vector $(2,0,-2)$.
\begin{enumerate}
\item If $\mathscr{E}$ is locally free, then $H^{i}(X,\mathscr{E})=0$ 
for $i=0,1,2$. \item If $\mathscr{E}$ is not locally free,
then $h^{1}(X,\mathscr{E})=h^{2}(X,\mathscr{E})\neq 0$.
\end{enumerate}
\end{prop}

\proof By Lemma \ref{lem:no0} we have $H^{0}(X,\mathscr{E})=0$ and
$h^{1}(X,\mathscr{E})=h^{2}(X,\mathscr{E})$. If $\mathscr{E}$ is locally free,
by Serre's duality $h^{2}(X,\mathscr{E}^{*})=0$. Then 
$h^{2}(X,\mathscr{E})=0$ by Lemma \ref{lem:edual}, and the first item is shown.
If $\mathscr{E}$ is not locally free, we have two cases.

\textit{Case 1}: $[\mathscr{E}]\in B\cap M_{10}^{s}$. Then
$\mathscr{E}^{**}=\mathscr{O}_{X}\oplus\mathscr{O}_{X}$, and we
have a short exact sequence
$$0\longrightarrow\mathscr{E}\longrightarrow
\mathscr{O}_{X}\oplus\mathscr{O}_{X}\longrightarrow\mathscr{G}\longrightarrow
0$$since $\mathscr{E}$ is torsion free, where $\mathscr{G}$ is
supported on a finite number of points. Thus
$h^{2}(X,\mathscr{E})=2$, and we are done.

\textit{Case 2}: $\mathscr{E}$ is strictly semistable. By Lemma 
1.1.5 in \cite{OG1}, $\mathscr{E}$ fits into an exact sequence
\begin{equation}
\label{eq:ezw}
0\longrightarrow \mathscr{I}_{Z}\longrightarrow\mathscr{E}
\longrightarrow\mathscr{I}_{W}\longrightarrow 0
\end{equation}
for some $Z,W\in Hilb^{2}(X)$. Since $H^{0}(X,\mathscr{I}_{W})=0$
and $h^{i}(X,\mathscr{I}_{Z})=1$ for $i=1,2$, the long exact
sequence induced by (\ref{eq:ezw}) implies
$h^{2}(X,\mathscr{E})\neq 0$.
\endproof

Let $\mathscr{F}$ be a universal family on $R\times X$, and 
consider the universal quotient module
\begin{equation}
\label{eq:es2}0\longrightarrow\mathscr{G}\longrightarrow
p_{X}^{*}\mathscr{H}\stackrel{\rho}\longrightarrow
\mathscr{F}\longrightarrow 0,
\end{equation}
where $p_{X}$ is the projection on $X$ and $\mathscr{H}:=H^{0}(X,\mathscr{E}(NH))
\otimes\mathscr{O}_{X}(-NH)$ for $N\in\mathbb{Z}$ sufficiently
big, where $\mathscr{E}$ is any sheaf parameterized by $M_{10}$.
In particular, $\mathscr{H}$ is locally free and
$H^{0}(X,\mathscr{H})=H^{1}(X,\mathscr{H})=0$. Notice
that any $s\in R$ corresponds to an exact sequence
\begin{equation}
\label{eq:es1}0\longrightarrow \mathscr{K}\longrightarrow
\mathscr{H}\stackrel{f_{\mathscr{E}}}\longrightarrow\mathscr{E}\longrightarrow
0.
\end{equation}
As $\mathscr{F}$ and $p_{X}^{*}\mathscr{H}$ are $R-$flat, the sheaf
$\mathscr{G}$ is $R-$flat. For any $s\in R$ let
$\mathscr{G}_{s}$ (resp. $(p_{X}^{*}\mathscr{H})_{s}$,
$\mathscr{F}_{s}$) denote the restriction of $\mathscr{G}$ (resp.
$p_{X}^{*}\mathscr{H}$, $\mathscr{F}$) to the fiber of the
projection $p_{R}:R\times X\longrightarrow R$ over the point $s$. Then
$$\mathscr{G}_{s}\simeq\ker((p_{X}^{*}\mathscr{H})_{s}\longrightarrow
\mathscr{F}_{s})=\ker(f_{\mathscr{E}})=\mathscr{K}.$$

\begin{prop}
\label{prop:ri}We have the following properties:
\begin{enumerate}
\item For every $i\in\mathbb{Z}$ the sheaves
$\mathbb{R}^{i}p_{R*}\mathscr{G}$ and
$\mathbb{R}^{i}p_{R*}(p_{X}^{*}\mathscr{H})$ are locally free of
rank $h^{i}(X,\mathscr{H})$. \item For every $s\in R$ and
$i\in\mathbb{Z}$, the canonical morphism
$$(\mathbb{R}^{i}p_{R*}\mathscr{F})_{s}\longrightarrow
H^{i}(p_{R}^{-1}(s),\mathscr{F}_{s})\simeq
H^{i}(X,\mathscr{E})$$is an isomorphism, where $\mathscr{E}$ is a
sheaf corresponding to the point $s\in R$.
\end{enumerate}
\end{prop}

\proof The main ingredient is the following lemma:

\begin{lem}
\label{lem:mum}Let $f:T\longrightarrow S$ be a proper morphism
of Noetherian schemes, and suppose $S$ to be reduced. Let
$\mathscr{U}\in Coh(T)$ be an $S-$flat family of sheaves, and let
$i\in\mathbb{Z}$. The function sending any $s\in S$ to
$h^{i}(T_{s},\mathscr{U}_{s})$ is constant if and only if 
$\mathbb{R}^{i}f_{*}\mathscr{U}$ is locally free and the
canonical morphism
$$(\mathbb{R}^{i}f_{*}\mathscr{U})_{s}\longrightarrow
H^{i}(T_{s},\mathscr{U}_{s})$$is an isomorphism. If this 
is verified, then $(\mathbb{R}^{i-1}f_{*}\mathscr{U})_{s}\longrightarrow
H^{i-1}(T_{s},\mathscr{U}_{s})$ is an isomorphism for every $s\in
S$.
\end{lem}

\proof See \cite{M}, Chapter II.5, Corollary 2.\endproof

We only need to show the proposition for $i=0,1,2$. For every $s\in R$ we
have $(p_{X}^{*}\mathscr{H})_{s}\simeq\mathscr{H}$, so that
$H^{i}(p_{R}^{-1}(s),(p_{X}^{*}\mathscr{H})_{s})\simeq
H^{i}(X,\mathscr{H})$, and the function sending $s\in R$ to
$h^{i}(p_{R}^{-1}(s),(p_{X}^{*}\mathscr{H})_{s})$ is constant. By
Lemma \ref{lem:mum}, the sheaf
$\mathbb{R}^{i}p_{R*}(p_{X}^{*}\mathscr{H})$ is locally free of
rank $h^{i}(X,\mathscr{H})$. In particular, as
$h^{0}(X,\mathscr{H})=h^{1}(X,\mathscr{H})=0$, then
$\mathbb{R}^{0}p_{R*}(p_{X}^{*}\mathscr{H})=\mathbb{R}^{1}
p_{R*}(p_{X}^{*}\mathscr{H})=0$.

The next step is to study $\mathbb{R}^{i}p_{R*}\mathscr{G}$.
Applying $\mathbb{R}p_{R*}$ to the exact sequence (\ref{eq:es2}), by the first
part of the proposition we get
$\mathbb{R}^{0}p_{R*}\mathscr{G}=0$ and
$\mathbb{R}^{1}p_{R*}\mathscr{G}\simeq\mathbb{R}^{0}p_{R*}\mathscr{F}$. 
We show that this last sheaf is trivial. Let $\mathscr{E}$ be 
a sheaf parameterized by $M_{10}$, and consider 
a corresponding point $s\in R$. Then $\mathscr{F}_{s}\simeq\mathscr{E}$, 
and the map sending $s$ to $H^{0}(X,\mathscr{F}_{s})$ is constant and 
trivial by Lemma \ref{lem:no0}. The canonical morphism
$(\mathbb{R}^{0}p_{R*}(\mathscr{F}))_{s}\longrightarrow
H^{0}(X,\mathscr{F}_{s})=0$ is then an isomorphism by Lemma
\ref{lem:mum}, so that $\mathbb{R}^{0}p_{R*}(\mathscr{F})=0$. It
remains to show that $\mathbb{R}^{2}p_{R*}\mathscr{G}$ is a
vector bundle of rank $h^{2}(X,\mathscr{H})$: consider $s\in R$ 
and its associated exact sequence (\ref{eq:es1}).
The long exact sequence induced by this and Lemma \ref{lem:no0}
imply $h^{2}(X,\mathscr{G}_{s})=h^{2} (X,\mathscr{H})$, so
that $\mathbb{R}^{2} p_{R*}\mathscr{G}$ is a vector bundle of rank
$h^{2}(X,\mathscr{H})$; for every $s\in R$ the canonical
morphism $(\mathbb{R}^{2}p_{R*}\mathscr{G})_{s}\longrightarrow
H^{2}(X,\mathscr{G}_{s})\simeq H^{2}(X,\mathscr{K})$ is an
isomorphism by Lemma \ref{lem:mum}.

Finally, we study $\mathbb{R}^{i}p_{R*}\mathscr{F}$ for $i=1,2$. As
$\mathbb{R}^{3}p_{R*}\mathscr{F}=0$, by Lemma \ref{lem:mum} the
canonical morphism
$(\mathbb{R}^{2}p_{R*}\mathscr{F})_{s}\longrightarrow
H^{2}(X,\mathscr{E})$ is an isomorphism. Let
$\xi:\mathbb{R}^{1}p_{R*}\mathscr{F}\longrightarrow
\mathbb{R}^{2}p_{R*}\mathscr{G}$ be the morphism induced by the
exact sequence (\ref{eq:es2}). Since
$\mathbb{R}^{1}p_{R*}(p_{X}^{*}\mathscr{H})=0$ by the first part
of the proof, $\xi$ is injective. In particular, for any $s\in R$
the morphism $\xi_{s}$ is injective, so that
$$(\mathbb{R}^{1}p_{R*}(\mathscr{F}))_{s}\simeq\ker(
(\mathbb{R}^{2}p_{R*}(\mathscr{G}))_{s}\stackrel{\delta}\longrightarrow
(\mathbb{R}^{2}p_{R*}(p_{X}^{*}\mathscr{H}))_{s}).$$The morphism 
$\delta$ is simply the morphism $H^{2}(X,\mathscr{K})\longrightarrow
H^{2}(X,\mathscr{H})$ induced by the exact sequence
(\ref{eq:es1}), by the previous part of the proof. Since
$H^{1}(X,\mathscr{H})=0$, we have $\ker(\delta)\simeq
H^{1}(X,\mathscr{E})$, and we are done.\endproof

We are finally able to prove the following

\begin{prop}
\label{prop:nol}We have $2B=\lambda(u_{2}(1))$.
\end{prop}

\proof By definition $\widetilde{\lambda}(u_{2}(1))=
det(\mathbb{R}p_{R*}(\mathscr{F}))$. By Theorem \ref{thm:2b}, the line bundle
$\widetilde{\lambda}(u_{2}(1))$ descends to $2B+\lambda(u_{1}(L))$ for some $L\in
Pic(X)$. Applying $\mathbb{R}p_{R*}$ to the exact
sequence (\ref{eq:es2}), by point 1 of Proposition \ref{prop:ri}
we get the exact sequence
$$0\longrightarrow\mathbb{R}^{1}p_{R*}(\mathscr{F})\longrightarrow
\mathbb{R}^{2}p_{R*}(\mathscr{G})\stackrel{\beta}\longrightarrow
\mathbb{R}^{2}p_{R*}(p_{X}^{*}\mathscr{H})\longrightarrow
\mathbb{R}^{2}p_{R*}(\mathscr{F})\longrightarrow 0.$$As
$\mathbb{R}^{0}p_{R*}\mathscr{F}=0$ by point 2 of Proposition
\ref{prop:ri}, we get $det(\mathbb{R}p_{R*}\mathscr{F})\simeq
det(\mathbb{R}^{2}p_{R*}(p_{X}^{*}\mathscr{H}))\otimes
det(\mathbb{R}^{2}p_{R*}\mathscr{G})^{-1}$. Then
$det(\beta)$ gives a
section $s$ of the line bundle $det(\mathbb{R}p_{R*}(\mathscr{F}))$, 
whose zero locus is given by the set where $det(\beta)$ is not an isomorphism. 
By Propositions \ref{prop:coholf} and \ref{prop:ri} this locus is exactly 
$p^{-1}(B)$, and we are done.\endproof

\subsection{Description of $H^{2}(\widetilde{M}_{10},\mathbb{Z})$}

Yoshioka (see Theorem 0.1 in \cite{Y}) showed the following: 
if $S$ is any projective K3 surface,  $v\in H^{2*}(S,\mathbb{Z})$ 
is a primitive Mukai vector with $(v,v)>0$ and $H$ is a $v-$generic polarization,
the moduli space $M_{v}$ of $H-$semistable sheaves on $S$ with Mukai vector 
$v$ is an irreducible symplectic variety, and there is an isometry of Hodge 
structures $v^{\perp}\longrightarrow H^{2}(M_{v},\mathbb{Z}),$ where 
$v^{\perp}$ is a sublattice of the Mukai lattice of $S$ and $H^{2}(M_{v},
\mathbb{Z})$ is a lattice with respect to the Beauville form.
In this section, we show an analogue of this in the case of
$\widetilde{M}_{10}$. Here, $X$ is a
projective K3 surface with Picard group spanned by an ample line bundle $H$
such that $H^{2}=2$, and $v=(2,0,-2)\in H^{2*}(X,\mathbb{Z})$.

\begin{lem}
\label{lem:h2}
We have $v^{\perp}\simeq H^{2}(X,\mathbb{Z})\oplus\mathbb{Z}$.
\end{lem}

\proof By the Hirzebruch-Riemann-Roch Theorem, $w\in H^{2*}(X,\mathbb{Z})$ is
orthogonal to $v$ if and only if $w=(r,c,r)$ for $r\in\mathbb{Z}$ 
and $c\in H^{2}(X,\mathbb{Z})$.
\endproof
 
We have the Hodge morphism $\widetilde{\mu}:
H^{2}(X,\mathbb{Z})\longrightarrow H^{2}(\widetilde{M}_{10},\mathbb{Z})$
respecting the lattice structures, and the morphism
$$c_{1}\circ\pi^{*}\circ\lambda\circ u:Pic(X)\oplus\mathbb{Z}\longrightarrow
H^{2}(\widetilde{M}_{10},\mathbb{Z}).$$By Proposition \ref{prop:eqdonpot}, 
these two morphisms agree on $Pic(X)$. Let
$$f:v^{\perp}\longrightarrow H^{2}(\widetilde{M}_{10},\mathbb{Z}),\,\,\,\,\,\,\,\
f(r,c,r):=\widetilde{\mu}(c)+c_{1}(\pi^{*}(\lambda(u_{2}(r)))).$$

\begin{thm}
\label{thm:isohodge}The morphism $f$ is a Hodge isometry between
$v^{\perp}$, viewed as a sublattice of the Mukai lattice
$\widetilde{H}(X,\mathbb{Z})$, and its image in
$H^{2}(\widetilde{M}_{10},\mathbb{Z})$, being a lattice with respect to the
Beauville form $q$.
\end{thm}

\proof The morphism $f$ is an injective
morphism of Hodge structures by point 1 of Theorem
\ref{thm:rapagnetta2}. By Proposition \ref{prop:nol}, $\lambda(u_{2}(r))=2rB$, so that
$$\pi^{*}\lambda(u_{2}(r))=2r\widetilde{B}+m\widetilde{\Sigma}$$
for some $m\in\mathbb{Z}$. Intersecting with $\delta$, by point 2 of Theorem \ref{thm:rapagnetta2} we
get $m=r$. In conclusion, we have $f(r,c,r)=\widetilde{\mu}(c)+2rc_{1}
(\widetilde{B})+rc_{1}(\widetilde{\Sigma}).$
By point 4 of Theorem \ref{thm:rapagnetta2} and by definition of 
the Mukai pairing, it is then an easy calculation to show that $f$
is an isometry.
\endproof

\section{The local factoriality of $M_{6}$}

From now on, we deal with the $6-$dimensional O'Grady's example
$M_{6}$, and we show that it is $2-$factorial. In this section we
recall the construction of $M_{6}$ and of $\widetilde{M}_{6}$, and
we resume the basic properties we need for the proof of the
$2-$factoriality. Moreover, we show that $M_{6}$ is not locally
factorial.

\subsection{Generalities on $M_{6}$}

In the following, let $C$ be a smooth projective curve of degree
$2$, and let $J:=Pic^{0}(C)$ be its jacobian surface. Suppose
there is an ample line bundle $H$ on $J$ such that $NS(J)=
\mathbb{Z}\cdot c_{1}(H)$ and $c_{1}^{2}(H)=2$. Finally, let
$\widehat{J}:=Pic^{0}(J)$ be the abelian surface dual to $J$.

Let $v:=(2,0,-2)\in\widetilde{H}(J,\mathbb{Z})$, and let $M_{v}$
be the moduli space of $H-$semistable sheaves on $J$ whose Mukai
vector is $v$. The regular locus of $M_{v}$ is the open subset $M_{v}^{s}$ 
parameterizing stable sheaves. Let $\Sigma_{v}$ be the singular locus
of $M_{v}$, which is a closed subset of codimension 2 in $M_{v}$
(see \cite{OG2}). Since in this setting any semistable locally
free sheaf is stable (see Lemma 2.1.2 in \cite{OG2}), the open subset $M_{v}^{lf}$
of $M_{v}$ parmaterizing locally free sheaves is contained in $M_{v}^{s}$. 
Let $\overline{B}_{v}$ be the closed subset of $M_{v}$ 
parameterizing non-locally free sheaves.
In particular, $\Sigma_{v}\subseteq \overline{B}_{v}$. Finally, let
$$a_{v}:M_{v}\longrightarrow J\times\widehat{J},\,\,\,\,\,\,\,\,\,\,\,
a_{v}([\mathscr{E}]):=\bigg(\sum
c_{2}(\mathscr{E}),det(\mathscr{E})\bigg),$$ and let
$M_{6}:=a_{v}^{-1}(0,\mathscr{O}_{J})$.

\begin{thm}(\textbf{O'Grady}, \textbf{Lehn-Sorger}). The moduli
space $M_{v}$ admits a symplectic resolution
$\pi_{v}:\widetilde{M}_{v} \longrightarrow M_{v},$
which is obtained as the blow-up of $M_{v}$ along $\Sigma_{v}$ with reduced 
schematic structure. Let $\widetilde{M}_{6}:=\pi_{v}^{-1}(M_{6})$. Then 
$\widetilde{M}_{6}$ is an irreducible symplectic variety 
of dimension 6 and second Betti number 8.
\end{thm}

\proof The proof is in \cite{OG2}. In \cite{L-S} it is show that
$\widetilde{M}_{v}$ can be obtained as the blow up of $M_{v}$
along its reduced singular locus.\endproof

Let $\widetilde{\Sigma}_{v}$ be the exceptional
divisor of $\pi_{v}$, and let $\widetilde{B}_{v}$ be the proper transform of
$\overline{B}_{v}$ under $\pi_{v}$. Let
$\pi:=\pi_{v|\widetilde{M}_{6}}$, and let $\Sigma:=\Sigma_{v}\cap
M_{6}$, the singular locus of $M_{6}$. In particular, $\pi$ is the
blow up of $M_{6}$ along $\Sigma$ with its reduced structure. Finally,
let $\overline{B}:=\overline{B}_{v}\cap M_{6}$,
$\widetilde{\Sigma}:=\pi^{-1}(\Sigma)$ (the exceptional divisor of
$\pi$) and $\widetilde{B}:=
\widetilde{B}_{v}\cap\widetilde{M}_{6}$ (the proper transform of
$\overline{B}$ under $\pi$). As shown in section 5.1 in
\cite{OG2}, $\widetilde{B}$ is an irreducible Weil divisor on
$\widetilde{M}_{6}$. Let
$M_{6}^{\mu-ss}$ be the Donaldson-Uhlenbeck compactification of the
moduli space $M_{6}^{\mu}$ of $\mu-$stable sheaves, and
let $\phi:M_{6}\longrightarrow M_{6}^{\mu-ss}$ be the canonical surjective
morphism. Let $\delta$ be the fiber of $\pi$ over
a generic point in the smooth locus of $\Sigma$, and let $\gamma$
be as in section 5.1 of \cite{OG2}. Finally, let $$\mu_{D}:H^{2}(J,\mathbb{Z})\longrightarrow
H^{2}(M_{6}^{\mu-ss},\mathbb{Z})$$be the Donaldson morphism. 

\begin{thm}
\label{thm:rapac}(\textbf{Rapagnetta}). Let $\widetilde{\mu}:=\pi^{*}
\circ\phi^{*}\circ\mu_{D}$.
\begin{enumerate}
\item The morphism
$\widetilde{\mu}:H^{2}(J,\mathbb{Z})
\longrightarrow H^{2}(\widetilde{M}_{6},\mathbb{Z})$ is injective.
\item There is a line bundle $A\in Pic(\widetilde{M}_{6})$ such
that $c_{1}(\widetilde{\Sigma})= 2c_{1}(A).$ \item We have the
following equalities:
$$c_{1}(A)\cdot\delta=-1,\,\,\,\,\,\,\,\,\,\,c_{1}(\widetilde{B})\cdot\delta=1,$$
$$c_{1}(A)\cdot\gamma=1,\,\,\,\,\,\,\,\,\,\,c_{1}(\widetilde{B})\cdot\gamma=-2.$$
\item The second integral cohomology of
$\widetilde{M}_{6}$ is
$$H^{2}(\widetilde{M}_{6},\mathbb{Z})=\widetilde{\mu}(H^{2}(J,\mathbb{Z}))
\oplus\mathbb{Z}\cdot c_{1}(A)\oplus\mathbb{Z}\cdot
c_{1}(\widetilde{B}).$$ \item Let $q$ be the Beauville form of $\widetilde{M}_{6}$.
Then for every $\alpha,\beta\in H^{2}(J,\mathbb{Z})$ we have
$$q(\widetilde{\mu}(\alpha),\widetilde{\mu}(\beta))=\alpha\cdot\beta,\,\,\,\,\,\,\,\,
q(\widetilde{\mu}(\alpha),c_{1}(\widetilde{B}))=q(\widetilde{\mu}
(\alpha),c_{1}(A))=0,$$
$$q(c_{1}(A),c_{1}(A))=-2,\,\,\,\,\,\,\,\,\,\,\,q(c_{1}(A),c_{1}(\widetilde{B}))=2,$$
$$q(c_{1}(\widetilde{B}),c_{1}(A))=2,\,\,\,\,\,\,\,\,\,\,\,
q(c_{1}(\widetilde{B}),c_{1}(\widetilde{B}))=-4.$$
\end{enumerate}
\end{thm}

\proof Item 1 is \cite{OG2}, Proposition 7.3.3. The proof of the
other points is contained in \cite{R2}, Theorems 3.3.1, 3.4.1 and
3.5.1.\endproof

\subsection{Flat families}

In this subsection we present two examples of flat families of sheaves
we will use in the following. As in section 2.2, we refer to section
10 for the general construction.

\begin{esem}
\label{esem:restj}Let $E$ be a rank 2 vector bundle on $J$ with
trivial first and second Chern classes and such that $hom(E,E)=2$.
Moreover, let $J[2]$ be the set of $2-$torsion points in $J$, and
let $y\in J\setminus J[2]$. Fix a surjective morphism
$\varphi:E\longrightarrow\mathbb{C}_{y}$ and let $\mathscr{K}:=\ker
(\varphi)$: by Lemma 4.3.3 in \cite{OG2}, any sheaf defining a
point in $\widetilde{B}_{v}$ is the kernel of a surjective
morphism from $\mathscr{K}$ to $\mathbb{C}_{x}$ for some point
$x\in J$. Let $p_{1},p_{2}:J\times J\longrightarrow J$ be the two
projections. As in Example \ref{esem:univx}, the sheaf
$p_{1*}\mathscr{H}om(p_{2}^{*}\mathscr{K},\mathscr{O}_{\Delta})$
is a vector bundle of rank 2, and for any $x\in J$ the canonical
morphism
$$(p_{1*}\mathscr{H}om(p_{2}^{*}\mathscr{K},\mathscr{O}_{\Delta}))_{x}
\longrightarrow Hom(\mathscr{K},\mathbb{C}_{x})$$is an isomorphism
(see Lemma \ref{lem:mum}). Let $Y:=\mathbb{P}(p_{1*}\mathscr{H}om(p_{2}^{*}\mathscr{K},
\mathscr{O}_{\Delta}))\stackrel{p}\longrightarrow J.$ There is a
tautological morphism (see Section 10)
$$\widetilde{f}:q_{J}^{*}\mathscr{K}\otimes q_{Y}^{*}\mathscr{T}
\longrightarrow (p\times id_{J})^{*}\mathscr{O}_{\Delta},$$whose
kernel is denoted $\mathscr{H}$.
\end{esem}

\begin{lem}
\label{lem:hj}Let $\mathscr{E}$ be a sheaf defining a point in
$\widetilde{B}_{v}$ whose bidual is $E$ and whose singular locus
is given by $x,y\in J$. Let
$f_{\mathscr{E}}:\mathscr{K}\longrightarrow\mathbb{C}_{x}$ be the
surjective morphism whose kernel is $\mathscr{E}$. Then
$f_{\mathscr{E}}$ defines a point $[f_{\mathscr{E}}]\in Y$, and
$\mathscr{H}_{[f_{\mathscr{E}}]}\simeq\mathscr{E}$. Moreover,
$\mathscr{H}$ is a $Y-$flat family and
$\widetilde{f}$ is surjective.
\end{lem}

\proof The proof is the same as the one of Lemma
\ref{lem:ssu}.\endproof

\begin{esem}
\label{esem:restjx}Let $E$ be as in the previous example, with the
further property that $det(E)\simeq\mathscr{O}_{J}$. Let $x\in J$
and $\varphi:E\longrightarrow \mathbb{C}_{-x},$ a surjective morphism
whose kernel is denoted $\mathscr{K}$. Let
$Y:=\mathbb{P}(p_{x*}\mathscr{H}om(p_{J}^{*}\mathscr{K},i_{*}\mathbb{C}_{x}))
\stackrel{p}\longrightarrow\{x\},$ where $p_{J}:\{x\}\times
J\longrightarrow J$ and $p_{x}:\{x\}\times J\longrightarrow\{x\}$
are the two projections, and $i:\{x\}\longrightarrow J$ is the
closed immersion. Then, $Y\simeq\mathbb{P}^{1}$, and its points
correspond to surjective morphisms from $\mathscr{K}$ to
$\mathbb{C}_{x}$. As before, we get a tautological
morphism
$\widetilde{f}:q_{J}^{*}\mathscr{K}\longrightarrow
j_{*}\mathscr{O}_{\mathbb{P}^{1}}(1),$ where
$j:\mathbb{P}^{1}\times\{x\}\longrightarrow\mathbb{P}^{1}\times J$
is the immersion. Let $\mathscr{H}:=\ker(\widetilde{f})$.
\end{esem}

\begin{lem}
\label{lem:hjx}Let $\mathscr{E}$ be a sheaf defining a point in
$\widetilde{B}$ whose bidual is $E$ and whose singular locus is
given by $x,-x\in J$. Let
$f_{\mathscr{E}}:\mathscr{K}\longrightarrow\mathbb{C}_{x}$ be the
surjective morphism whose kernel is $\mathscr{E}$. Then
$f_{\mathscr{E}}$ defines a point $[f_{\mathscr{E}}]\in Y$, and
$\mathscr{H}_{[f_{\mathscr{E}}]}\simeq\mathscr{E}$. Moreover,
$\mathscr{H}$ is a $Y-$flat family and $\widetilde{f}$ is surjective.
\end{lem}

\proof Again, the proof is the same as the one of Lemma
\ref{lem:ssu}, using the Claim in section 5.1 of \cite{OG2}.
\endproof

\subsection{The moduli space $M_{6}$ is not locally factorial}

A first application of Theorem \ref{thm:rapac} is the following:

\begin{lem}
\label{lem:d}There is a non-trivial irreducible Weil divisor $D\in
A^{1}(M_{6})$ such that $2D=0$. If $\widetilde{D}$ is the proper
transform of $D$ by $\pi$, then there is $m\in\mathbb{Z}$ such
that $A=\widetilde{D}+m\widetilde{\Sigma}$ in the group
$Div(\widetilde{M}_{6})$ of Weil divisors of $\widetilde{M}_{6}$.
\end{lem}

\proof As in the proof of Proposition \ref{prop:picm}, we have
$A^{1}(M_{6})\simeq Pic(\pi^{-1}(M_{6}^{s}))$. The restriction
of $A$ to $\pi^{-1}(M_{6}^{s})$ defines then an irreducible Weil
divisor $D\in A^{1}(M_{6})$. By point 2 of Theorem
\ref{thm:rapac} we have
$$2D=2A_{|\pi^{-1}(M_{6}^{s})}=\widetilde{\Sigma}_{|\pi^{-1}(M_{6}^{s})}=0.$$Now,
the Weil divisor $\widetilde{\Sigma}$ is a prime divisor, so it is
a generator for the group $Div(\widetilde{M}_{6})$. Since $A$ is a
line bundle on $\widetilde{M}_{6}$, it defines an element in
$Div(\widetilde{M}_{6})$, so that there are
$m,m_{1},...,m_{n}\in\mathbb{Z}$ and prime divisors
$D_{1},...,D_{n}$ such that
$$A=m\widetilde{\Sigma}+\sum_{i=1}^{n}m_{i}D_{i}.$$As
$A_{|\pi^{-1}(M_{6}^{s})}=\sum_{i=1}^{n}m_{i}D_{i|\pi^{-1}(M_{6}^{s})}$,
we have $\sum_{i=1}^{n}m_{i}D_{i}=\widetilde{D}$, and we are done.
It remains to show that $D$ is not trivial: if $D=0$, then
$\widetilde{D}=0$, so that $A=m\widetilde{\Sigma}=2mA$ (by point 2
of Theorem \ref{thm:rapac}). Then $c_{1}(A)$ is torsion in
$H^{2}(\widetilde{M}_{6},\mathbb{Z})$, which is not possible 
by point 4 of Theorem \ref{thm:rapac}.
\endproof

\begin{prop}
\label{prop:m6nf}The Weil divisor $D$ is not Cartier, and $M_{6}$
is not locally factorial.
\end{prop}

\proof If $D$ was a Cartier divisor, then $\pi^{*}(D)=
\widetilde{D}+kA,$ for some $k\in\mathbb{Z}$. By Lemma \ref{lem:d}
we then get $\pi^{*}(D)=(1-2m+k)A.$ The integer $1-2m+k$ is odd:
indeed, if there was $n\in\mathbb{Z}$ such that $2n=1-2m+k$, then
$\pi^{*}(D)=n\widetilde{\Sigma}$ and we would have
$$D=\pi^{*}(D)_{|\pi^{-1}(M_{6}^{s})}=n\widetilde{\Sigma}_{|\pi^{-1}
(M_{6}^{s})}=0,$$which is not possible since $D$ is non-trivial.
By point 3 of Theorem
\ref{thm:rapac} and the fact that $\delta$ is contracted by $\pi$,
one gets
$$0=c_{1}(\pi^{*}(D))\cdot\delta=(1-2m+k)c_{1}(A)\cdot\delta=2m-k-1.$$As
$2m-k-1\neq 0$, we get a contradiction, and $D$ is not a
Cartier divisor. Finally, this clearly implies that $M_{6}$ cannot
be locally factorial.
\endproof

\begin{oss}
\label{oss:notor6} As a consequence of this, $Pic(M_{6})$ is free. 
Indeed, let $L\in Pic(M_{6})$ be torsion of period $t$, and let
$\widetilde{L}$ be its proper transform under $\pi$. Then
$\pi^{*}(L)=\widetilde{L}+kA$ for some $k\in\mathbb{Z}$, and
$t(\widetilde{L}+kA)=0$. As $Pic(\widetilde{M}_{6})$ has no
torsion by point 4 of Theorem \ref{thm:rapac}, we get
$\widetilde{L}=-kA$, and
$$L=\widetilde{L}_{|\pi^{-1}(M_{6}^{s})}= -kA_{|\pi^{-1}(M_{6}^{s})}=
-kD.$$As $L\in Pic(M_{6})$, we get $kD\in Pic(M_{6})$, so that $k$
has to be even by Proposition \ref{prop:m6nf} and Lemma
\ref{lem:d}. In conclusion $L=0$, and we are done.

The same proof even shows that $\pi^{*}:Pic(M_{6})\longrightarrow
Pic(\widetilde{M}_{6})$ is injective. As in Remark
\ref{oss:notor10}, from this one can deduce that the morphism
$c_{1}:Pic(M_{6})\longrightarrow H^{2}(M_{6},\mathbb{Z})$ is
injective.
\end{oss}

\section{Line bundles on $M_{6}$}

In this section we calculate the Picard groups of
$\widetilde{M}_{6}$ and of $M_{6}$ following the same argument
as in sections 3 and 4.

\subsection{Le Potier's construction}

Let $e:=[\mathscr{E}]\in K_{top}(J)$ be the class of
a sheaf $\mathscr{E}$ parameterized by $M_{6}$, and let $h:=[H]\in
K_{top}(J)$.

\begin{lem}
\label{lem:kj}Let $\alpha\in K_{top}(J)$. Then $\alpha\in e^{\perp}\cap
\{1,h,h^{2}\}^{\perp\perp}$ if and only if
$c_{1}(\alpha)\in c_{1}(H)^{\perp\perp}$ and $ch_{2}(\alpha)=rk(\alpha)
\eta_{J}\in H^{4}(J,\mathbb{Z})$, where $\eta_{J}$ is the fundamental class of $J$.
\end{lem}

\proof The proof works as the one of Lemma \ref{lem:perp}.\endproof

Using this lemma, we are able to prove the following:

\begin{prop}
\label{prop:lineb6}Let $p\in J$ be any point, and let $$u_{1}:Pic(J)\longrightarrow e^{\perp}\cap
\{1,h,h^{2}\}^{\perp\perp},\,\,\,\,\,\,\,\,u_{1}(L):=
[\mathscr{O}_{J}-L]+\frac{c_{1}^{2}(L)}{2}[\mathbb{C}_{p}],$$
$$u_{2}:\mathbb{Z}\longrightarrow
e^{\perp}\cap\{1,h,h^{2}\}^{\perp\perp},\,\,\,\,\,\,\,\,\, u_{2}(n):=
n[\mathscr{O}_{J}]+n[\mathbb{C}_{p}].$$The morphism $u:=u_{1}+u_{2}$ 
is a group isomorphism.
\end{prop}

\proof The proof works as the one of Proposition \ref{prop:lineb}.
\endproof

\begin{prop}
\label{prop:lep6}Let $i:M_{6}\longrightarrow M_{v}$ be the
inclusion. There is a group morphism $\widetilde{\lambda}_{6}:=i^{*}\circ\lambda\circ
u:Pic(J)\oplus\mathbb{Z}\longrightarrow Pic(M_{6}),$ where $\lambda$ is the Le
Potier morphism. In particular, this induces a group morphism
$$\lambda_{6}:NS(J)\oplus\mathbb{Z}\longrightarrow Pic(M_{6})$$such that
for any $L\in Pic(J)$, $n\in\mathbb{Z}$ we have $\lambda_{6}(c_{1}(L),n)=
\widetilde{\lambda}_{6}(L,n)$.
\end{prop}

\proof The existence of the maps $\lambda\circ u$ and
$\widetilde{\lambda}_{6}$ is implied by Lemma \ref{lem:kj} and
Theorem 8.1.5 in \cite{H-L}. The fact that $\lambda\circ u$, and
hence $\widetilde{\lambda}_{6}$, is a group morphism is as in the proof
of Proposition \ref{prop:lineb}. As $v:K_{top}(J)\longrightarrow 
H^{2*}(J,\mathbb{Z})$ is an isomorphism, if $c_{1}(L_{1})=c_{1}(L_{2})$, 
then $u(L_{1})=u(L_{2})$, and we are done. But this implies the existence 
of the morphism $\lambda_{6}$ defined on $NS(J)\oplus\mathbb{Z}$, and we are done.
\endproof

In the following, let $\widetilde{\lambda}_{6,1}:=i^{*}\circ\lambda\circ
u_{1}$ and $\widetilde{\lambda}_{6,2}:=i^{*}\circ\lambda\circ u_{2}$,
so that $\widetilde{\lambda}_{6}=\widetilde{\lambda}_{6,1}+
\widetilde{\lambda}_{6,2}$. Then $\widetilde{\lambda}_{6,1}$ induces a
morphism $\lambda_{6,1}:NS(J)\longrightarrow Pic(M_{6})$, such that for 
every $L\in Pic(J)$ we have $\lambda_{6,1}(c_{1}(L))=\widetilde{\lambda}_{6,1}(L)$.
Then we have $\lambda_{6}=\lambda_{6,1}+\widetilde{\lambda}_{6,2}$.

\begin{lem}
\label{lem:l6nogamma}We have the following intersection properties.
\begin{enumerate}
\item Let $L\in Pic(J)$. Then
$c_{1}(\pi^{*}\widetilde{\lambda}_{6,1}(L)))\cdot\gamma=
c_{1}(\pi^{*}\widetilde{\lambda}_{6,1}(L))))\cdot\delta=0.$
\item Let $n\in\mathbb{Z}$. Then 
$c_{1}(\pi^{*}\widetilde{\lambda}_{6,2}(n))\cdot\gamma=-n$.
\end{enumerate}
\end{lem}

\proof We begin with the first item.
The equality $c_{1}(\pi^{*}\widetilde{\lambda}_{6}(L))\cdot
\delta=0$ is trivial, as $\delta$ is contracted by $\pi$. Notice
that
$$c_{1}(\pi^{*}\widetilde{\lambda}_{6,1}(L))\cdot\gamma=
c_{1}(\pi^{*}\widetilde{\lambda}_{6,1}(L)_{|\gamma})=c_{1}
(\lambda_{\mathscr{H}}(u_{1}(L)))$$by Theorem 8.1.5 in \cite{H-L} and
Lemma \ref{lem:hjx}, where $\lambda_{\mathscr{H}}$ is the Le
Potier's morphism defined using the flat family $\mathscr{H}$ of
Example \ref{esem:restjx}. By the Grothendieck-Riemann-Roch
Theorem, we have
$$c_{1}(\lambda_{\mathscr{H}}(u_{1}(L)))=q_{Y*}[q_{J}^{*}(ch(u_{1}(L))td(J))
\cdot ch(\mathscr{H}))]_{3}\in
H^{2}(\mathbb{P}^{1},\mathbb{Z})=$$ $$=-q_{Y*}
(q_{J}^{*}(c_{1}(L))\cdot ch_{2}(\mathscr{H}))=0,$$and we are done.
For the second item, the proof is the same as the one of point 2 of 
Proposition \ref{prop:noint}, using Example \ref{esem:restjx}.\endproof

\subsection{Donaldson's and Le Potier's morphisms}

The main result we need is the following:

\begin{prop}
\label{prop:l6don}For any $L\in Pic(J)$ we have
$c_{1}(\pi^{*}\widetilde{\lambda}_{6,1}(L))=\widetilde{\mu}(c_{1}(L))$.
\end{prop}

\proof The proof of this proposition is almost the same as the one
of Proposition \ref{prop:eqdonpot}. Let $L\in Pic(J)$ and let $Y$
and $\mathscr{H}$ be as in Example \ref{esem:restj}. Using the
same argument as in Step 1 of the proof of Proposition
\ref{prop:eqdonpot}, we get
\begin{equation}
\label{eq:uguali}c_{1}(\pi^{*}\lambda(u_{1}(L)))_{|Y}=
\widetilde{\mu}(c_{1}(L))_{|Y}.
\end{equation}
Now, let $Y_{6}:=Y\cap\widetilde{M}_{6}$, and let
$\beta:=c_{1}(\pi^{*}\widetilde{\lambda}_{6,1}(L))-
\widetilde{\mu}(c_{1}(L))\in
H^{2}(\widetilde{M}_{6},\mathbb{Z}).$ By equation
(\ref{eq:uguali}), we have $\beta_{|Y_{6}}=0$, and by point 1 of
Lemma \ref{lem:l6nogamma} and the definition of $\widetilde{\mu}$ we
have $\beta\cdot\gamma=\beta\cdot\delta=0$. Following Step 2 of
the proof of Proposition \ref{prop:eqdonpot}, these two properties
imply $\beta=0$, and we are done.\endproof

\begin{cor}
\label{cor:picmt6}The morphism $\lambda_{6}:NS(J)\oplus\mathbb{Z}
\longrightarrow Pic(M_{6})$ is injective. Moreover, we have
$Pic(\widetilde{M}_{6})=\pi^{*}\lambda_{6,1}(NS(J))\oplus\mathbb{Z}\cdot[A]
\oplus\mathbb{Z}\cdot[\widetilde{B}].$
\end{cor}

\proof The proof works as that of Corollary \ref{cor:pic},
using Proposition \ref{prop:l6don}.\endproof

\section{The $2-$factoriality of $M_{6}$}

We are now able to show the $2-$factoriality of $M_{6}$. We
need to add a remark on $\overline{B}$: the proper transform of
$\overline{B}$ is an irreducible Weil divisor in
$\widetilde{M}_{6}$, so that
$\overline{B}=\Sigma\cup B$ for some irreducible Weil divisor $B$
of $M_{6}$ whose proper transform is $\widetilde{B}$.

\begin{prop}
\label{prop:pim6}We have $A^{1}(M_{6})=\lambda_{6}(NS(J))\oplus\mathbb{Z}\cdot[B]\oplus
\mathbb{Z}/2\mathbb{Z}\cdot[D]$.
\end{prop}

\proof The proof is similar to the one of Proposition \ref{prop:picm}, and
we need to show that the following sequence
$$0\longrightarrow\mathbb{Z}\stackrel{\sigma}\longrightarrow
Pic(\widetilde{M}_{6})\stackrel{\rho}\longrightarrow
Pic(\pi^{-1}(M_{6}^{s}))\longrightarrow 0$$is exact, where $\sigma(1):=\widetilde{\Sigma}$
and $\rho$ is the restriction morphism. The only thing to prove is that if $L\in
\ker(\rho)$, then it is a multiple of
$\widetilde{\Sigma}$. By Corollary \ref{cor:picmt6}, there are
$M\in Pic(J)$ and $n,m\in\mathbb{Z}$ such that
$L=\pi^{*}(\lambda_{6,1}(c_{1}(M)))+n\widetilde{B}+mA.$ By Lemma
\ref{lem:d}, we have $\rho(L)=\lambda_{6,1}(c_{1}(M))+nB+mD\in A^{1}(M_{6}).$ As
$\rho(L)=0$, then $\rho(2L)=0$,
so that $2nB=\lambda_{6,1}(2c_{1}(M))$, as $2mD=0$ by Lemma
\ref{lem:d}. In particular, their proper transforms are equal,
getting $2n\widetilde{B}=\pi^{*}(\lambda_{6,1}(2c_{1}(M)))$, so that
$$-4n=2nc_{1}(\widetilde{B})\cdot\gamma=c_{1}(\pi^{*}\lambda_{6,1}(2c_{1}(M))\cdot\gamma
=0,$$by point 3 of Theorem \ref{thm:rapac} and point 1 of Lemma
\ref{lem:l6nogamma}. In conclusion, $n=0$ and
$\lambda_{6,1}(2c_{1}(M))=0$. By Corollary \ref{cor:picmt6} then $c_{1}(M)=0$. 
In conclusion $L=mA$ for some $m\in\mathbb{Z}$, so that
$$0=\rho(L)=\rho(mA)=mD.$$By Lemma
\ref{lem:d}, then, $m$ is even and $L$ is a multiple of
$\widetilde{\Sigma}$.
\endproof

Here is the main result of this section:

\begin{thm}
\label{thm:m6lf}There is a line bundle $L\in Pic(J)$ and $t\in
\mathbb{Z}/2\mathbb{Z}$ such that $B+tD=\lambda_{6}(c_{1}(L),1).$ 
In particular, $M_{6}$ is $2-$factorial.
\end{thm}

\proof By Lemma \ref{lem:kj} and Proposition \ref{prop:pim6} there
are $M\in Pic(J)$, $n\in\mathbb{Z}$ and
$t\in\mathbb{Z}/2\mathbb{Z}$ such that
$\widetilde{\lambda}_{6,2}(1)=\lambda_{6,1}(c_{1}(M))+nB+tD\in
A^{1}(M_{6}).$ In particular $nB+tD\in Pic(M_{6})$: 
we need to show that $n=1$. Taking
the pull-back of $nB+tD$ to $\widetilde{M}_{6}$ there is
$m\in\mathbb{Z}$ such that
$$n\widetilde{B}+mA=\pi^{*}(nB+tD)=\pi^{*}(\lambda_{6}(-c_{1}(M),1)).$$By 
point 3 of Theorem \ref{thm:rapac} we get
$$0=\pi^{*}(nB+tD)\cdot\delta=n\widetilde{B}\cdot\delta+mA\cdot\delta=n-m,$$as
$\delta$ is contracted by $\pi$, and
$$-2n+m=n\widetilde{B}\cdot\gamma+mA\cdot\gamma=\pi^{*}(\lambda_{6}(-c_{1}(M),1))
\cdot\gamma=-1$$by Lemma \ref{lem:l6nogamma}. In
conclusion, $n=1$ and we are done. It remains to show that 
$M_{6}$ is $2-$factorial: since $B+tD$ is a Cartier divisor, we have
$$\lambda_{6}(NS(J))\oplus\mathbb{Z}[B+tD]\subseteq
Pic(M_{6}).$$We have then
two possibilities: the first one is $t=0$, so that $B$ is
Cartier. In this case, the only Weil divisor which is not Cartier
is $D$, and we are done. The second case is $t=1$, so that $B+D$ is Cartier. 
As $2D=0$, we then get $2B\in Pic(M_{6})$, and we are done.
\endproof

\begin{oss}
\label{oss:pbbd}As seen in the proof, one has
$\pi^{*}(\lambda_{6,2}(1))=\widetilde{B}+A+
\pi^{*}\lambda_{6,2}(c_{1}(L))$ for some line bundle $L\in Pic(J)$. 
As it was pointed out to me by Rapagnetta, using our construction one
can easily show that there is a line bundle $A\in
Pic(\widetilde{M}_{6})$ such that $2A=\widetilde{\Sigma}$. Indeed,
as shown in \cite{OG2}, we have
$$H^{2}(\widetilde{M}_{6},\mathbb{Q})=\widetilde{\mu}(H^{2}(J,\mathbb{Q}))
\oplus\mathbb{Q}\cdot c_{1}(\widetilde{B})\oplus\mathbb{Q}\cdot
c_{1}(\widetilde{\Sigma}),$$so that there are $\beta\in
H^{2}(J,\mathbb{Q})$ and $p,q\in\mathbb{Q}$ such that
$$c_{1}(\pi^{*}\lambda_{6,2}(1))=
\widetilde{\mu}(\beta)+p\widetilde{B}+q\widetilde{\Sigma}.$$ By
equation 7.3.5 in \cite{OG2} one gets
$$0=c_{1}(\pi^{*}\lambda_{6,2}(1))\cdot\delta=p-2q,$$
$$-1=c_{1}(\pi^{*}\lambda_{6,2}(1))\cdot\gamma=
-2p+2q.$$In conclusion $q=1/2$ and $p=1$. Now,
$c_{1}(\pi^{*}\lambda_{6,2}(1))\in
H^{2}(\widetilde{M}_{6},\mathbb{Z})$, so that if
$\widetilde{\Sigma}$ was a generator for
$H^{2}(\widetilde{M}_{6},\mathbb{Z})$, we would have
$q\in\mathbb{Z}$, which is clearly not the case. Then, there must
be a line bundle $A\in Pic(\widetilde{M}_{6})$ such that
$2c_{1}(A)=c_{1}(\widetilde{\Sigma})$, and we are done.
\end{oss}

\section{The Beauville form of $\widetilde{M}_{6}$}

In this last section, we prove an analogue of Theorem
\ref{thm:isohodge} about the Beauville form of
$\widetilde{M}_{6}$. Here is the result:

\begin{thm}
\label{thm:isoh6}Let $v=(2,0,-2)\in\widetilde{H}(J,\mathbb{Z})$.
There is a morphism of Hodge structures
$$f:v^{\perp}\longrightarrow H^{2}
(\widetilde{M}_{6},\mathbb{Z}),$$which is an isometry between
$v^{\perp}$, as a sublattice of the Mukai lattice
$\widetilde{H}(J,\mathbb{Z})$, and its image in
$H^{2}(\widetilde{M}_{6},\mathbb{Z})$, lattice with respect to the
Beauville form $q$.
\end{thm}

\proof As in Lemma \ref{lem:h2}, a Mukai vector $w$ is
orthogonal to $v$ if and only if $w=(r,c,r)$ for $r\in\mathbb{Z}$
and $c\in H^{2}(J,\mathbb{Z})$, so that $v^{\perp}\simeq
H^{2}(J,\mathbb{Z})\oplus\mathbb{Z}$. Let
$$f:v^{\perp}\longrightarrow
H^{2}(\widetilde{M}_{6},\mathbb{Z}),\,\,\,\,\,\,\,\,\,\,f((r,c,r))
:=\widetilde{\mu}(c)+rc_{1}(\widetilde{B})+ rc_{1}(A).$$The
morphism $f$ is an injective morphism of Hodge structures.
By point 5 of Theorem \ref{thm:rapagnetta2} and 
definition of the Mukai pairing, it is easy to see that $f$
is an isometry on its image.
\endproof

\section{Appendix: Construction of flat families}

In this Appendix we resume a general construction of
flat families we used in several occasions.
Let $S$ be an algebraic surface and $T$ a proper scheme. Let
$p_{S}$ and $p_{T}$ be the two obvious projections of $T\times S$.
Let $\mathscr{V}$ and $\mathscr{W}$ be two $T-$flat coherent
sheaves on $T\times S$, and suppose
$p_{T*}\mathscr{H}om(\mathscr{V},\mathscr{W})$ to be a vector bundle
on $T$. Let
$Y:=\mathbb{P}(p_{T*}\mathscr{H}om(\mathscr{V},\mathscr{W}))$ and
$p:Y\longrightarrow T$ the projection morphism. Let $\mathscr{T}$ 
be the tautological line bundle on $Y$. As shown
in \cite{H}, Chapter II, Prop. 7.11, there is a canonical
morphism $f:\mathscr{T}\longrightarrow p^{*}
p_{T*}\mathscr{H}om(\mathscr{V},\mathscr{W})$ which is injective.
Finally, let $q_{Y}:Y\times S\longrightarrow Y$ and $q_{S}:Y\times
S\longrightarrow S$ be the two projections. We have the following
commutative diagram:

\begin{center}
$\begin{CD}Y @<{q_{Y}}<< Y\times S
@>{q_{S}}>> S\\@V{p}VV @VV{p\times id_{S}}V @|\\
T @<<p_{T}< T\times S @>>p_{S}> S
\end{CD}$
\end{center}

\noindent and the following equality holds:
$$p^{*}p_{T*}\mathscr{H}om(\mathscr{V},\mathscr{W})=q_{Y*}\mathscr{H}om((p\times
id_{S})^{*}\mathscr{V},(p\times id_{S})^{*}\mathscr{W}).$$By
the projection formula $f$ defines then a global section $$\sigma\in
H^{0}(Y\times S,\mathscr{H}om((p\times
id_{S})^{*}\mathscr{V}\otimes q_{Y}^{*}\mathscr{T},(p\times
id_{S})^{*}\mathscr{W})),$$corresponding to a morphism
$$\widetilde{f}:(p\times id_{S})^{*}\mathscr{V}\otimes
q_{Y}^{*}\mathscr{T}\longrightarrow (p\times
id_{S})^{*}\mathscr{W}.$$This construction allows us to produce
flat families, as shown in section 2.2 and 6.2.

\subsection*{Acknowledgment}
This work is the issue of my PhD Thesis at the Universit\'e de Nantes,
under the direction of Christoph Sorger. I would like to thank
him for his support, useful discussions and improvements. I also would
like to thank Antonio Rapagnetta and Baohua Fu for having helped me in
several occasions, showing me different approaches and improvements.
Finally, I would like to thank Manfred Lehn, Eyal Markman and Dimitri 
Markushevich for having read my PhD thesis and having signaled me
mistakes and improvements.

Arvid Perego, Institut f\"ur Mathematik der Universit\"at Johannes Gutenberg Mainz,
Staudinger Weg 9, 55099 Mainz;
E-mail address: perego@uni-mainz.de 

\end{document}